\documentclass[12pt]{elsart}
\usepackage{amsfonts}
\usepackage{amsmath}
\usepackage{graphicx}
\usepackage{epstopdf}

\newcommand{\R}{\Bbb{R}}

\newcommand{\N}{\Bbb{N}}
\newcommand{\sign}{\mbox{\rm sign}}

\newtheorem{theorem}{Theorem}[section]
\newtheorem{proposition}[theorem]{Proposition}
\newtheorem{definition}[theorem]{Definition}
 \newtheorem{lemma}[theorem]{Lemma}
\newtheorem{corollary}[theorem]{Corollary}

\newtheorem{remark}[theorem]{Remark}

\journal{}

\begin{document}

\begin{frontmatter}
\title{Global behaviour of a second order nonlinear difference equation}
\author{Ignacio Bajo} 
and
%\ead{ibajo@dma.uvigo.es}
\author{Eduardo Liz
%\thanksref{dges}
}
%\footnote{Author for correspondence}}
%\ead{eliz@dma.uvigo.es}
%\ead[url]{http://www.dma.uvigo.es/\~{} eliz/}
%\fntext[corresp]{Author for correspondence} 
\address{Departamento de Matem\'atica Aplicada II, E.T.S.E. Telecomunicaci\'on,
Universidade de
Vigo, Campus Marcosende, 36310 Vigo, Spain}
%\thanks[dges]{Author for correspondence, eliz@dma.uvigo.es.}

\begin{abstract}
\noindent 
We describe the asymptotic behaviour and the stability properties of the solutions to the nonlinear second order difference equation 
$$
x_{n+1}=\frac{x_{n-1}}{a+bx_nx_{n-1}},\, n\geq 0,
$$
for all values of the real parameters $a,b$, and any initial condition $(x_{-1},x_0)\in\R^2$.
\end{abstract}

\end{frontmatter}

%%%
%%%
%%%
\section{Introduction}
%%%
%%%
%%%

We consider the second order difference equation defined by
\begin{equation}
\label{1}
x_{n+1}=\frac{x_{n-1}}{a+bx_nx_{n-1}},
\end{equation}
for all integers $n\geq 0$, with initial condition $(x_{-1},x_0)\in\R^2$, where $a,b$ are real parameters. 

This recurrence has recently attracted some attention, and several particular cases were already studied. Notice that the difference equation
$$
x_{n+1}=\frac{\delta x_{n-1}}{\beta+\gamma x_nx_{n-1}}
$$
can be obviously reduced to (\ref{1}) if $\delta\neq 0$.

As far as we know, the first particular case of (\ref{1}) considered in the literature was $a=b=1$, with positive initial conditions. {\c{C}}inar \cite{ci1}  established a formula for the subsequences of even and odd terms of the solutions, respectively. 
Stevi{\'c} \cite{st} gave further insight for this case, showing that every solution of (\ref{1}) converges to zero if $x_{-1}x_0\neq 0$ and $x_{-1}x_0\neq -1/n$, for all positive integer $n$. If $x_{-1}x_0= 0$, then Eq. (\ref{1}) is $2$-periodic.

The case when  $a>0$, $b>0$, and the initial conditions are nonnegative was also considered by  {\c{C}}inar  \cite{ci2,ci4}, who stated a similar formula for the subsequences of even and odd terms of the solutions. Later,  Andruch-Sobilo and Migda \cite{am} proved that these subsequences are convergent; moreover, it is shown that  they converge to zero if $a\geq 1$.

The case $a=-1$, $b>0$, and arbitrary initial conditions such that $bx_{-1}x_0\neq 1$ was addressed in \cite{ci3,ci5}. The author finds the representation formula as in the previous cases and proves that for $bx_{-1}x_0>1$, one of the subsequences converges to zero and the other one diverges.

Aloqeili \cite{al} investigated (\ref{1}) in the case when $a>0$, $b=-1$, and proved some interesting results on the asymptotic behaviour of the solutions. He shows that they typically converge to zero if $a>1$ and are oscillating if $0<a<1$. 

Andruch-Sobilo and Migda \cite{am2} considered the case $a<0$, $b>0$, with nonnegative initial conditions, showing that  the subsequences of even and odd terms are monotone.

In this paper, we succeed in giving a complete picture of the asymptotic behaviour of the solutions of (\ref{1}) depending on the involved parameters and the initial data. 

An advantage of our approach is that it allows us to consider arbitrary initial conditions, and to split our study in three cases, depending only on the values of the parameter $a$. We get as particular cases all the mentioned known results in the literature about Equation (\ref{1}), but also other cases are solved here for the first time.  Our results on bifurcation and stability of the solutions  in sections 4 and 5 are also new.

The key feature of the solutions to (\ref{1}) is that the sequence $\{y_n\}$ defined by $y_n=x_nx_{n-1}$ for all $n\geq 0$ solves a rational difference equation of M\"obius type
\begin{equation}
\label{2}
y_{n+1}=\frac{y_{n}}{a+by_n},
\end{equation}
which is reducible to a linear difference equation (see, e.g., \cite{cull,el}).

We notice that other authors refer to Equation (\ref{2}) as a  {\it Ricatti difference equation} \cite[Section 3.3]{ag}, \cite[Section 1.6]{kl}. When $a>0, b>0$, Equation (\ref{2})  is equivalent to the Pielou logistic equation \cite[Example 3.3.4]{ag}, \cite[Example 2.39]{el}.

This fact allows us to write Equation (\ref{1}) in the form 
$$
x_{n+1}=h(n)x_{n-1},
$$
where the term $h(n)$ depends only on the parameter $a$ and the product $\alpha=bx_{-1}x_0$ for each $n\geq 0$. Thus, the subsequences of even and odd terms from a solution $\{x_n\}$ of (\ref{1}) are given by the expressions
$$
x_{2k+2}=x_0\prod_{i=0}^kh(2i+1)\quad ;\quad x_{2k+1}=x_{-1}\prod_{i=0}^kh(2i),\; k\geq 0.
$$

Notice that if these subsequences converge, say, $\displaystyle\lim_{k\to\infty}x_{2k+1}=p$, $\displaystyle\lim_{k\to\infty}x_{2k+2}=q$, then, by continuity arguments, the pair $(p,q)$ satisfies
\begin{equation}
\label{pq}
q=\frac{q}{a+bqp}\quad ;\quad p=\frac{p}{a+bqp}\, .
\end{equation}
Hence, either $p=q=0$ or $pq=(1-a)/b$. In particular, $\{p,q,p,q,\dots\}$ is a $2$-periodic solution of (\ref{1}), in such a way that the solution $\{x_n\}$  either converges to zero or to a $2$-periodic solution. For this reason, the analysis of the convergence of $\{x_{2k}\}$ and  $\{x_{2k+1}\}$ is an important step in our proofs. 

The paper is organized as follows: in Section 2 we derive the mentioned representation for the solutions of (\ref{1}). In Section 3 we describe the asymptotic behaviour of the solutions; it is divided into three subsections depending on the values of $a$.  In Section 4 we give an interpretation of our results in terms of a bifurcation problem. Finally, we devote Section 5 to analyze the stability properties of the periodic solutions of (\ref{1}). 

%%%
%%%
%%%
\section{A formula for the solutions}
%%%
%%%
%%%

Throughout the paper, we denote $\alpha=bx_{-1}x_0$. In this section, we state a representation formula for the solutions of (\ref{1}) starting at any initial condition $(x_{-1},x_0)\in\R^2$, except the following cases:
\begin{enumerate}
\item $a=1$ and $\alpha= -1/n$,  for some $n\geq 1$.
\item  $a\neq 1$ and $\alpha= a^{n}(a-1)/(1-a^{n})$, for some $n\geq 1$.
\end{enumerate}

We emphasize that, in these cases, it is not possible to construct a complete solution $\{x_n\}_{n={-1}}^{\infty}$ starting at $(x_{1},x_0)$, since at some point the denominator in (\ref{1}) becomes zero.

On the other hand, it is convenient to consider the cases $\alpha=0$ and $\alpha=1-a$ separately due to their singularity (see Propositions \ref{1p} and \ref{2p} below).

These facts motivate us to introduce the following definitions:

\begin{definition} 
\label{1d}
We say that the pair $(x_{-1},x_0)$ is an admissible initial condition for (\ref{1}) if either $a=1$ and $\alpha\neq -1/n$, or $a\neq 1$ and $\alpha\neq a^{n}(a-1)/(1-a^{n})$, for all $n\geq 1$.

Solutions of (\ref{1}) corresponding to admissible initial conditions are called admissible solutions.
\end{definition}

\begin{definition} 
An admissible solution of (\ref{1}) is called a regular solution if $\alpha\neq 0$ and $\alpha\neq 1-a$. Admissible solutions that are not regular are called singular solutions.
\end{definition}

In the following two propositions, we describe the singular solutions of (\ref{1}). We notice that for $\alpha=0$ all solutions are admissible if $a\neq 0$, while for $\alpha=1-a$ all solutions are admissible.

\begin{proposition}
\label{1p}
Assume that $\alpha=0$ and $a\neq 0$. Then, $x_{2k}=x_0/a^{k}$ and $x_{2k-1}=x_{-1}/a^{k}$, for all $k\geq 1$.
\end{proposition}
{\bf Proof.}
First, assume that $b=0$. Hence,  Eq. (\ref{1}) reduces to $x_{n+1}=x_{n-1}/a$, and it follows by induction that $x_{2k}=x_0/a^{k}$ and $x_{2k-1}=x_{-1}/a^{k}$, for all $k\geq 1$.

If $x_0=0$, then 
$$
x_2=\frac{x_{0}}{a+bx_1x_{0}}=0.
$$
It is easily derived by induction that $x_{2k}=0$ for all $k\geq 1$.

On the other hand,
$$
x_{2k+1}=\frac{x_{2k-1}}{a+bx_{2k}x_{2k-1}}=\frac{x_{2k-1}}{a}, \, \forall \, k\geq 0.
$$
Thus, $x_{2k+1}=x_{-1}/a^{k+1}$, for all $k\geq 0$. 

The proof in the case $x_{-1}=0$ is completely analogous. In this case we obtain $x_{2k-1}=0$, $x_{2k}=x_{0}/a^k$, for all $k\geq 1$. 
\qed\\

\begin{proposition}
\label{2p}
If $\alpha=1-a$, then the solution of (\ref{1}) with initial condition $(x_{-1},x_0)$ is $2$-periodic.
\end{proposition}
{\bf Proof.} Using (\ref{1}), we have
\begin{align*}
x_1&=\frac{x_{-1}}{a+bx_0x_{-1}}=\frac{x_{-1}}{a+\alpha}=x_{-1},\\
\noalign{\smallskip}
x_2&=\frac{x_{0}}{a+bx_1x_{0}}=\frac{x_{0}}{a+bx_{-1}x_{0}}=\frac{x_{0}}{a+\alpha}=x_0.
\end{align*}
The result follows by induction.
\qed\\

In order to get a representation for the regular solutions of (\ref{1}), we first  state the explicit expression for the solutions of the M\"obius recurrence (\ref{2}). For a proof of the following proposition, see, e.g., \cite[Example 2.39]{el}.
\begin{proposition}
\label{p2} Assume that $\alpha=by_0$ is in the conditions of Definition \ref{1d}, and $\alpha\not\in\{ 0,1-a\}$.
Then, the solution of Equation  (\ref{2}) starting at the initial condition $y_0$ is given by
$$
y_n=\left\{
\begin{array}{cl} \displaystyle\frac{y_0(1-a)}{a^n(1-a)+y_0b(1-a^n)} & \mbox{if $a\neq 1$;}\\
\noalign{\medskip}\displaystyle\frac{y_0}{1+y_0bn} & \mbox{if $a= 1$.}\end{array}
  \right.
$$
\end{proposition}

Now we are in a position to provide a representation for all admissible solutions of (\ref{1}).
\begin{theorem}
\label{t1}
Denote $\alpha=bx_{-1}x_0$. 
If $(x_{-1},x_0)$ is an admissible initial condition for (\ref{1}), then
the corresponding solution is given by
\begin{align*}
x_{2k+2}&=x_0\prod_{i=0}^kh(2i+1)\\
x_{2k+1}&=x_{-1}\prod_{i=0}^kh(2i),
\end{align*}
for all integer $k\geq 0$, 
where
\begin{equation}
\label{hn}
h(n)=\left\{
\begin{array}{cl} \displaystyle\frac{a^n(1-a)+\alpha(1-a^n)}{a^{n+1}(1-a)+\alpha(1-a^{n+1})} & \mbox{if $a\neq 1$;}\\
\noalign{\medskip}\displaystyle\frac{1+\alpha n}{1+\alpha (n+1)} & \mbox{if $a= 1$.}
\end{array}
  \right.
\end{equation}
\end{theorem}

The proof of Theorem \ref{t1}  follows by induction from the following result:

\begin{proposition}
\label{co1}
 If $\{x_n\}$ is an admissible solution of (\ref{1}), then 
\begin{equation}
\label{fh}
x_{n+1}=h(n)x_{n-1},
\end{equation}
for all $n\geq 0$, where $h(n)$ is defined by (\ref{hn}).
\end{proposition}
{\bf Proof.} First, we assume that $\{x_n\}$ is a regular admissible solution of (\ref{1}).

Denote $y_n=x_nx_{n-1}$. Multiplying Equation (\ref{1}) in both sides by $x_n$, we get
$$
y_{n+1}=x_{n+1}x_n=\frac{x_nx_{n-1}}{a+bx_nx_{n-1}}=\frac{y_{n}}{a+by_n}.
$$
On the other hand,
$$
x_{n+1}=\frac{x_{n-1}}{a+bx_nx_{n-1}}=\frac{x_{n-1}}{a+by_n}=h(n)x_{n-1},
$$
where $h(n)=1/(a+by_n)$. 
Since $\{y_n\}$ is the solution of (\ref{2}) with initial condition $y_0=x_{-1}x_0$,
a direct application of Proposition \ref{p2} gives formula (\ref{hn}).

For singular solutions, the result is also true. From Proposition \ref{1p}, we can see that  (\ref{fh}) holds for $\alpha=0$ with $h(n)=1/a$ for all $n\geq 0$. Using Proposition \ref{2p}, it is clear that (\ref{fh}) is satisfied for  $\alpha=1-a$ with $h(n)=1$ for all $n\geq 0$.
\qed\\

\begin{remark}
The formulas given for some particular cases of Equation (\ref{1}) in references  \cite{al,am,ci1,ci2,ci3,ci4,ci5,st} are particular cases of Theorem \ref{t1}. 
\end{remark}

%%%
%%%
%%%
\section{Asymptotic behaviour of the solutions}
%%%
%%%
%%%

In this section, we use the representation formula given in Section 2 to study the asymptotic behaviour of the solutions to (\ref{1}). We will consider three different cases.

%%%
\subsection{The case $a=-1$}
%%%

When $a=-1$, the expression for all admissible solutions given in Theorem \ref{t1} becomes very simple:
\begin{align}
\label{m1}
x_{2k+2}&=x_0(\alpha-1)^{k+1}\\
\label{m12}
x_{2k+1}&=\displaystyle\frac{x_{-1}}{(\alpha-1)^{k+1}},
\end{align}
for all $k\geq 0$. We notice that all solutions with $\alpha\neq 1$ are admissible. Moreover, they are regular if $\alpha\not\in\{0,2\}$.

Thus, we have the following result:
\begin{theorem}
\label{tm1}
If $a=-1$, then all regular solutions of (\ref{1}) are unbounded. Moreover, if $\{x_n\}$ is a regular solution of (\ref{1}), then some subsequences of $\{x_n\}$ are divergent and other converge to zero. The singular solutions are $2$-periodic if $\alpha=2$, and $4$-periodic if $\alpha=0$.
\end{theorem}
{\bf Proof.} 
It follows easily from the relations (\ref{m1})-(\ref{m12}). The complete behaviour of the solutions is the following (here $\sign(x)=1$ if $x>0$, and $\sign(x)=-1$ if $x<0$):
\begin{enumerate}
\item If $\alpha>2$ then $\displaystyle\lim_{k\to\infty}x_{2k+2}=\sign(x_0)\, \infty$, $\displaystyle\lim_{k\to\infty}x_{2k+1}=0.$
\item If $\alpha=2$ then $x_{2k+2}=x_0$,  $x_{2k+1}=x_{-1},$ for all $k\geq 0$.
\item If $1<\alpha<2$ then $\displaystyle\lim_{k\to\infty}x_{2k+2}=0$, $\displaystyle\lim_{k\to\infty}x_{2k+1}=\sign(x_{-1})\, \infty.$
\item If $0<\alpha<1$ then $\displaystyle\lim_{k\to\infty}x_{2k+2}=0$, $\displaystyle\lim_{k\to\infty}x_{4k+1}=-\sign(x_{-1})\, \infty,$ and $\displaystyle\lim_{k\to\infty}x_{4k+3}=\sign(x_{-1})\, \infty.$
\item If $\alpha=0$ then $x_{4k+1}=-x_{-1}$,  $x_{4k+2}=-x_0$, $x_{4k+3}=x_{-1}$,   $x_{4k+4}=x_{0},$ for all $k\geq 0$. Thus, $x_{k+4}=x_k$, for all $k\geq -1$.
\item If $\alpha<0$ then $\displaystyle\lim_{k\to\infty}x_{4k}=\sign(x_0) \infty$, $\displaystyle\lim_{k\to\infty}x_{4k+2}=-\sign(x_0) \infty,$ and  $\displaystyle\lim_{k\to\infty}x_{2k+1}=0.$
\end{enumerate}
\qed

%%%
\subsection{The case $|a|\geq1$, $a\neq -1$}
%%%

The main result in this section is the following:
\begin{theorem}
\label{tam1}
If $|a|\geq 1$ and $a\neq -1$, then all admissible solutions of (\ref{1}) converge to zero if  $\alpha\neq 1-a$, and are $2$-periodic if $\alpha=1-a$.
\end{theorem}
{\bf Proof.} 
 From Proposition \ref{2p}, we already know that  the solutions of (\ref{1}) are $2$-periodic if $\alpha=1-a$.
 Next we assume that $\alpha\neq 1-a$.

We first address the case $a=1$.
As we already mentioned in the introduction, it was proved in \cite{st} that all regular solutions of (\ref{1}) converge to zero if $a=b=1$. The same arguments of Theorem 1 in \cite{st} apply to the case $a=1$, $b\neq 0$, using Theorem \ref{t1}.

It remains the case $|a|>1$, $\alpha\neq 1-a$.
Let $\{x_n\}$ be  an admissible solution of (\ref{1}). Using Proposition \ref{co1}, and the fact that $\displaystyle\lim_{n\to\infty}a^{-n}=0$, we have:
\begin{align*}
\lim_{n\to\infty}\frac{|x_{n+1}|}{|x_{n-1}|}&=\lim_{n\to\infty} |h(n)|=\lim_{n\to\infty}\left|\frac{a^n(1-a)+\alpha(1-a^n)}{a^{n+1}(1-a)+\alpha(1-a^{n+1})}\right|\\
\noalign{\medskip}
&=\lim_{n\to\infty}\left|\frac{a^n(1-a-\alpha)+\alpha}{a^{n+1}(1-a-\alpha)+\alpha}\right|=
\lim_{n\to\infty}\left|\frac{(1-a-\alpha)+\alpha a^{-n}}{a(1-a-\alpha)+\alpha a^{-n}}\right|=\\\noalign{\medskip}
&=\frac{1}{|a|}<1.\end{align*}
The D'Alembert criterion ensures that $\displaystyle\lim_{k\to\infty}|x_{2k}|=\displaystyle\lim_{k\to\infty}|x_{2k+1}|=0$. Hence, $\displaystyle\lim_{n\to\infty}x_n=0$.
\qed

%%%
\subsection{The case $|a|<1$}
%%%

We begin this subsection with a simple result corresponding to the case $a=0$. Notice that in this case all solutions are admissible if $\alpha\neq 0$.
\begin{proposition}
\label{p0}
If $a=0$, then all admissible solutions of (\ref{1}) are $2$-periodic.
\end{proposition}
{\bf Proof.} 
In this case, Eq. (\ref{1}) becomes $x_{n+1}=1/(bx_n).$ Thus,
$x_{n+2}=1/(bx_{n+1})=x_n,$ for all $n\geq -1$. \qed\\

The proof of the following lemma is very easy from expression  (\ref{hn}), so we omit it:
\begin{lemma}
\label{lemah}
Let $a\neq 1$ and $\alpha\neq a^{n}(a-1)/(1-a^{n})$, for all $n\geq 1$.
Then, $h(n)=1-g(n)$, where
\begin{equation}
\label{gn}
g(n)=\frac{(a+\alpha-1)(1-a)a^n}{a^{n+1}(1-a)+\alpha(1-a^{n+1})},
\end{equation}
for all $n\geq 0$.
\end{lemma}

In order to address the case $0<|a|<1$, we investigate the character of the subsequences of even and odd terms, which depend on the sequence $\{h(n)\}$. For example, if $h(2k)>1$ for all sufficiently large $k$, then it is clear from (\ref{fh}) that the subsequence of odd terms is eventually  increasing.

Notice that, in view of Lemma \ref{lemah}, $h(n)<1$ if and only if $g(n)>0$, and  $h(n)>1$ if and only if $g(n)<0.$

\begin{proposition}
\label{ph}
Assume that $0<|a|<1$, $\alpha\neq 0$, $\alpha\neq 1-a$, and $\alpha\neq a^{n}(a-1)/(1-a^{n})$, for all $n\geq 1$. Then, there exists $N\in\N$ such that the sequences $\{g(2k)\}$ and $\{g(2k+1)\}$ have constant sign for all $k\geq N$.
\end{proposition}
{\bf Proof.} 
We first consider the case $a\in (0,1)$, and distinguish three situations:
\begin{enumerate}
\item If $a\in (0,1)$ and $\alpha>1-a$, then $(a+\alpha-1)(1-a)a^n>0$ and 
$a^{n+1}(1-a)+\alpha(1-a^{n+1})>(1-a)>0,$ for all $n\geq 0$. Thus, $g(n)>0$ for all $n\in\N$.
\item  If $a\in (0,1)$ and $0<\alpha<1-a$, then $(a+\alpha-1)(1-a)a^n<0$ and 
$a^{n+1}(1-a)+\alpha(1-a^{n+1})>\alpha>0,$ for all $n\geq 0$. Thus, $g(n)<0$ for all $n\in\N$.
\item  If $a\in (0,1)$ and $\alpha<0<1-a$, then $(a+\alpha-1)(1-a)a^n<0$ for all $n\geq 0$.  On the other hand, since  
$$\lim_{n\to\infty}a^{n+1}(1-a)+\alpha(1-a^{n+1})=\alpha<0,$$
 we can conclude that there exists  $n_0\in\N$ such that  $g(n)>0$ for all $n\geq n_0$.
\end{enumerate}
 Analogously, we consider the same situations for $a\in (-1,0)$.
\begin{enumerate}
\item If $a\in(-1,0)$ and $\alpha>1-a$, then $(a+\alpha-1)(1-a)a^{2k}>0$ and $(a+\alpha-1)(1-a)a^{2k+1}<0,$ for all $k\geq 0$. 
Since $a^{n+1}(1-a)+\alpha(1-a^{n+1})>(1-a)>0,$ for all $n\geq 0$,
 it follows that   $g(2k)>0$ and  $g(2k+1)<0$ for all $k\geq 0$.
\item  If $a\in (-1,0)$ and $0<\alpha<1-a$, then $(a+\alpha-1)(1-a)a^{2k}<0$ and $(a+\alpha-1)(1-a)a^{2k+1}>0,$ for all $k\geq 0$. Since
$$\lim_{n\to\infty}a^{n+1}(1-a)+\alpha(1-a^{n+1})=\alpha>0,
$$ it follows that there exists  $k_1\in\N$ such that  $g(2k)<0$ and  $g(2k+1)>0$ for all $k\geq k_1$.
\item  If $a\in (-1,0)$ and $\alpha<0<1-a$, then $(a+\alpha-1)(1-a)a^{2k}<0$ and $(a+\alpha-1)(1-a)a^{2k+1}>0,$ for all $k\geq 0$. Since
$$\lim_{n\to\infty}a^{n+1}(1-a)+\alpha(1-a^{n+1})=\alpha<0,
$$ it follows that there exists  $k_2\in\N$ such that  $g(2k)>0$ and  $g(2k+1)<0$ for all $k\geq k_2$.
\end{enumerate}
\qed

As a consequence of Proposition \ref{ph}, we have: 

\begin{corollary}
\label{cormon}
If $0<|a|<1$, and $\{x_n\}$ is a regular solution of (\ref{1}), then the subsequences $\{x_{2k}\}$ and $\{x_{2k+1}\}$ are eventually monotone.
\end{corollary}
{\bf Proof.} From Proposition  \ref{ph}, we know that  there exists $N\in\N$ such that the sequences $\{g(2k)\}$ and $\{g(2k+1)\}$ have constant sign for all $k\geq N$. Assume that $g(2k)>0$  for all $k\geq N$. Then, $h(2k)=1-g(2k)<1$  for all $k\geq N$. On the other hand, since $\lim_{n\to\infty}h(n)=1>0$, it is clear that there exists $N_1\geq N$ such that $0<h(2k)<1$ for all $k\geq N_1$.

Since, by  Proposition \ref{co1}, $x_{2k+1}=h(2k)x_{2k-1},\,\forall\, k\geq 0$, it follows that $x_{2k+1}<x_{2k-1},\,\forall\, k\geq N_1$, that is, the subsequence $\{x_{2k+1}\}_{k=N_1}^{\infty}$ is decreasing.

The remainder cases are analogous.
\qed\\

Using this corollary, we can prove the following key result:
\begin{proposition}
\label{pconv}
If $0<|a|<1$, and $\{x_n\}$ is a  regular solution of (\ref{1}), then the subsequences $\{x_{2k}\}$ and $\{x_{2k+1}\}$ are convergent.
\end{proposition}
{\bf Proof.} 
We only prove this result for the sequences of even terms, since the other case is completely analogous. 

As we noticed above, $\lim_{n\to\infty}h(n)=1>0$, and therefore $h(n)>0$ for all sufficiently large $n$. Without loss of generality, we assume that $h(n)>0$ for all $n\geq 0$.

Since $\{x_{2k}\}$ is eventually monotone, we only have to prove that it is bounded.

 For it, we use Theorem \ref{t1} and Lemma \ref{lemah}:
\begin{align*}
|x_{2k}|&=|x_0|\prod_{i=0}^{k-1}h(2i+1)=|x_0|\exp\left(\sum_{i=0}^{k-1}\ln(h(2i+1))\right)
\\
\noalign{\medskip}
&= |x_0|\exp\left(\sum_{i=0}^{k-1}\ln\left(1-\frac{(a+\alpha-1)(1-a)a^{2i+1}}{a^{2i+2}(1-a)+\alpha(1-a^{2i+2})}\right)\right)
\\
\noalign{\medskip}
&\leq |x_0|\exp\left((1-\alpha-a)(1-a)\sum_{i=0}^{k-1}\frac{a^{2i+1}}{a^{2i+2}(1-a)+\alpha(1-a^{2i+2})}\right)\\
\noalign{\medskip}
&:= |x_0|\exp\left((1-\alpha-a)(1-a)\sum_{i=0}^{k-1}\beta(i)\right).
\end{align*}
For the inequality above, we have used that $\ln(1-x)\leq -x$ for all $x<1$.

Since
\begin{align*}
\lim_{i\to\infty}\frac{|\beta(i+1)|}{|\beta(i)|}&=\lim_{i\to\infty}\left|\frac{a^{2i+2}\left(a^{2i+2}(1-a)+\alpha (1-a^{2i+2})\right)}{a^{2i+1}\left(a^{2i+3}(1-a)+\alpha (1-a^{2i+3})\right)} \right|\\
\noalign{\medskip}
&=\lim_{i\to\infty}\left|\frac{a\left(a^{2i+2}(1-a)+\alpha (1-a^{2i+2})\right)}{\left(a^{2i+3}(1-a)+\alpha (1-a^{2i+3})\right)} \right|=|a|<1,
\end{align*}
it follows from the D'Alembert rule that the series $\sum_{i=0}^{\infty}\beta(i)$ is convergent. This ensures that $\{x_{2k}\}$ is bounded.
\qed\\

Finally, we can state the main result of this subsection for the regular solutions of (\ref{1}). 

\begin{theorem}
\label{tap1}
If $0<|a|<1$, then all regular solutions of (\ref{1}) converge to a $2$-periodic solution $(p,q)$ of (\ref{1}), with $pq=(1-a)/b\neq 0$.
\end{theorem}
{\bf Proof.} 
By Proposition \ref{pconv}, there exist $\lim_{k\to\infty}x_{2k+1}=p\in\R$, $\lim_{k\to\infty}x_{2k+2}=q\in\R$. As it was mentioned in the introduction, the sequence $\{y_n\}=\{x_nx_{n-1}\}$ is a solution of (\ref{2}) and, by Proposition \ref{p2}, 
$$pq=\lim_{n\to\infty}y_n=\frac{1-a}{b}\neq 0.$$
Taking limits as $n\to\infty$ in (\ref{1}), it is clear that  the relations (\ref{pq}) hold, and therefore $(p,q)$ is a $2$-periodic solution of (\ref{1}). \qed\\

Using Theorem \ref{tap1} and Propositions \ref{1p}, \ref{2p} and \ref{p0}, we can describe completely the behaviour of all admissible solutions in the case $|a|<1$.

\begin{theorem}
\label{tap2}
Assume that $|a|<1$, and $\{x_n\}$ is an admissible solution of (\ref{1}). Then:
\begin{enumerate}
\item If either $a=0$ or $\alpha= 1-a$, then $\{x_n\}$ is $2$-periodic.
\item If $a\neq 0$, $\alpha\neq 0$ and $\alpha\neq 1-a$, then $\{x_n\}$ converges to a $2$-periodic solution.
\item If $a\neq 0$, $\alpha=0$ and $(x_{-1},x_0)\neq (0,0)$, then $\{x_n\}$ is unbounded.
\item If   $x_{-1}=x_0=0,$ then $x_n=0$ for all $n\geq 1$.
\end{enumerate}
\end{theorem}

\begin{remark}
\label{r311} Notice that zero is the unique equilibrium of (\ref{1}) if $(1-a)b\leq0$. When $(1-a)b>0$, there are two nontrivial equilibrium points $x_{\pm}=\pm((1-a)/b)^{1/2}$. Thus, the minimal period of the $2$-periodic solution $(p,q)$ mentioned in Theorem \ref{tap2} is actually $1$ if $(1-a)b>0$ and $p=q=x_{\pm}$.
\end{remark}

As a by-product of Theorems \ref{tm1}, \ref{tam1} and \ref{tap2}, we have the following result on the boundedness of the solutions to (\ref{1}):
\begin{proposition}
\label{pbound}
All admissible solutions of (\ref{1}) are bounded, except in the following two cases:
\begin{enumerate}
\item $a=-1$ and $\alpha\not\in\{ 0,2\}$;
\item $|a|<1$, $\alpha=0$, and $(x_{-1},x_0)\neq (0,0)$.
\end{enumerate}
\end{proposition}

%%%
\section{A bifurcation point of view}
%%%

The analysis made in Section 3 can be also viewed in terms of  bifurcation diagrams. First, notice that  the case $b>0$ may be reduced to  $b=1$ by the change of variables $v_n=b^{1/2}x_n$, and the case $b<0$ may be reduced to  $b=-1$ by the change of variables $v_n=(-b)^{1/2}x_n$. Thus, we can view Equation (\ref{1}) as a one-parameter family of difference equations depending only on $a$, if we consider the cases $b>0$, $b=0$, and $b<0$ separately.

As an example, we consider the case $b=-1$.

There are two {\it regular} bifurcation points in $a=-1$ and $a=1$. When $a>1$, all regular solutions converge to zero; as $a$ passes through $1$ to the left, the $\omega$-limit set of any regular solution is a $2$-periodic point. One of the branches of this periodic solution in the bifurcation diagram approaches zero as $a$ tends to $-1$, and the other one diverges to $+\infty$ or $-\infty$. After crossing the other bifurcation point $a=-1$, only the bounded branch remains, and all regular solutions are attracted by zero.

\begin{figure}[htp]
\centering
\includegraphics[totalheight=2in]{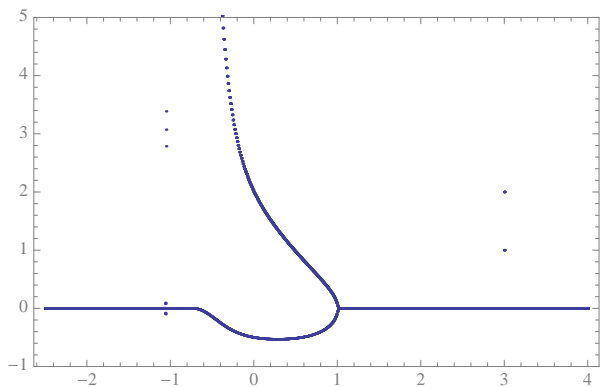}
\caption{Bifurcation diagram for (\ref{1}) with $b=-1$, $x_{-1}=1$, $x_0=2$.}
\label{bif1} 
\end{figure}

If we plot the bifurcation diagram corresponding to an initial condition $(x_{-1},x_0)$ with $\alpha\neq 0$, we also observe a {\it singular} bifurcation point when $1-a=\alpha=-x_{-1}x_0$ (that is, for $a^*=1+x_{-1}x_0$) if $|1+x_{-1}x_0|>1$. Indeed, for all values of $a$ in a neighbourhood of $a^*$ the limit of the solution is zero, while for $a=a^*$ the solution is $2$-periodic.
In figure \ref{bif1}, we plotted the bifurcation diagram corresponding to  $b=-1$ and the initial condition $(x_{-1},x_0)=(1,2)$. 
We observe the singular bifurcation point $a^*=3$, for which the solution is $2$-periodic.
\\

\begin{figure}[htp]
\centering
\includegraphics[totalheight=1.6in]{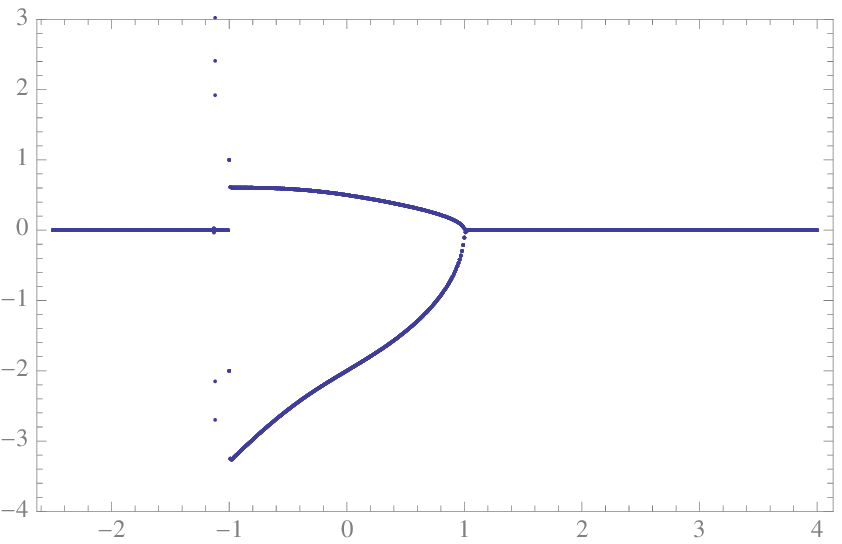}\hspace*{1cm}
\includegraphics[totalheight=1.6in]{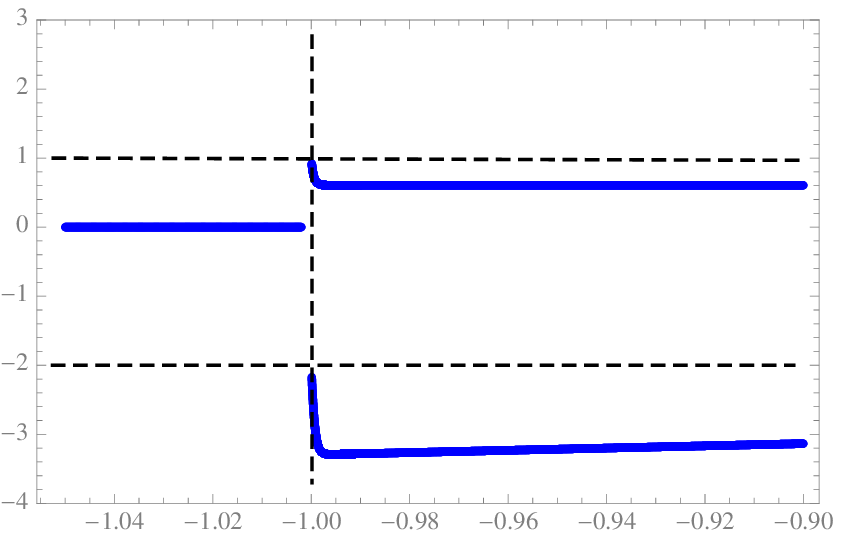}
\caption{Bifurcation diagrams for (\ref{1}) with $b=-1$, $x_{-1}=1$, $x_0=-2$.}
\label{bif2} 
\end{figure}

As it may be seen from Proposition \ref{pbound}, the case when $\alpha=2$ is special because all admissible solutions of (\ref{1}) are bounded. For $b=-1$, this happens when $x_{-1}x_0=-2$. We plot in Figure \ref{bif2} the bifurcation diagram corresponding to the initial data $(x_{-1},x_0)=(1, -2)$. On the right we plotted a magnification for $a$ close to $-1$ in order to emphasize that the branches of periodic points for this initial condition are continuous on the right at $a=-1$. Of course, there is a discontinuity on the left, since for $a=-1$ the solution is $2$-periodic, and for $a<-1$ it converges to zero.

\begin{remark} We produced the bifurcation diagrams in a standard way: for all values of $a$ with a step $0.005$, we constructed the first $400$ iterations corresponding to the initial conditions  $(x_{-1},x_0)=(1,2)$ and  $(x_{-1},x_0)=(1, -2)$, and plotted them starting at $n=350$. For the magnification in Figure \ref{bif2}, we produced $1600$ iterations for each value of $a$, with a step $0.0001$.
\end{remark}

%%%
%%%
%%%
\section{Stability properties}
%%%
%%%
%%%

As we have shown, the $\omega$-limit set of a bounded solution of (\ref{1}) is a periodic solution of minimal period $1$, $2$ or $4$. In this section, we study the stability properties of these periodic solutions. We begin with the zero solution.

\begin{proposition}
The zero solution of (\ref{1}) is asymptotically stable if and only if either $|a|>1$ or $a=1$ and $b\neq 0$. Moreover, in both cases  it attracts all regular solutions.
\end{proposition}
{\bf Proof.} 
The characteristic equation associated to the linearization of (\ref{1}) at the equilibrium $x=0$ is given by the quadratic equation $x^2=1/a$. Thus, the zero solution is locally asymptotically stable if $|a|>1$, and unstable if $|a|<1$. Theorem \ref{tam1} shows that zero is actually a global attractor of all regular solutions when $a\not\in [-1,1)$.

If $a=1$ and $b=0$, then all solutions are $2$-periodic (the minimal period may be one), and they are clearly stable but not asymptotically stable.

Finally, if $a=-1$, it follows from Theorem \ref{tm1} that zero is unstable, since solutions starting at initial conditions arbitrarily close to $(0,0)$ are unbounded.\qed\\

Next we deal with the nontrivial periodic solutions.

The unique periodic solutions with period greater than $2$ are the $4$-periodic points indicated in Theorem \ref{tm1} for $a=-1$ and $\alpha=0$. They are clearly unstable.

As proved in Proposition \ref{p0}, all admissible solutions of (\ref{1}) for $a=0$ are $2$-periodic; moreover, from the proof of this proposition, it is clear that they are stable. 

If $a\neq 0$, then the $2$-periodic points of (\ref{1}) are defined by the initial conditions $(p,q)$ such that $a+bpq=1$. This can be easily seen taking into account the correspondence between the solutions of (\ref{1})  and
the orbits of the discrete dynamical system associated to the map $F$ defined by 
$F(x,y)=(y, x/(a+bxy)).$ The $2$-periodic solutions of (\ref{1}) are defined by the fixed points of the map $F^2=F\circ F$. It is straightforward to prove that, if $a\neq 0$ and $(a,b)\neq (1,0)$,  $F^2(x,y)=(x,y)$ if and only if $a+bxy=1$. Notice that, as mentioned in Remark \ref{r311}, Equation (\ref{1}) has  two nontrivial equilibria $x_{\pm}=\pm((1-a)/b)^{1/2}$ if  $(1-a)b> 0$; thus, the minimal period of $(p,q)$ is one if $p=q=\pm((1-a)/b)^{1/2}$. Otherwise, the minimal period is two.

Direct computations show that the linearization of $F^2$ at any point $(p,q)$ satisfying $a+bpq=1$ has two eigenvalues: $a$ and $1$. Thus, all nonzero $2$-periodic solutions are unstable for $|a|>1$. When $a=-1$ or $a=1$, it follows from Theorems \ref{tm1}, \ref{tam1} that they are also unstable.

Since we already studied the case $a=0$, the remainder part of this section is devoted to prove that every nonzero periodic solution of (\ref{1}) is stable if $0<|a|<1$. For it, we will use the formula given in Theorem \ref{t1}. First, we need some bounds for the involved products.

\begin{lemma}
\label{ln1}
Assume that $0<|a|<1$ and $\alpha >(1-a)/2$. For all $n\geq 0$ one has:
$$a^n(1-a)+\alpha (1-a^n)>0,$$
and, as a consequence, $h(n)>0\, ,\, \forall\ n\geq 0.$
\end{lemma}
{\bf Proof.} 
Since $\alpha >(1-a)/2$, we have
$$a^n(1-a)+\alpha (1-a^n)>a^n(1-a)+(1-a^n)(1-a)/2=(1+a^n)(1-a)/2>0.\qed$$
%\qed\\

\begin{proposition}
\label{pn1} If $0<|a|<1$ and $\alpha >(1-a)/2$, then
\begin{enumerate}
\item[{\rm 1.}] $\displaystyle\prod_{i=0}^k h(2i+1)\leq \exp\left\{\displaystyle\frac{2|(1-a-\alpha)a|}{1-a^2}\right\}.$
\item[{\rm 2.}] $\displaystyle\prod_{i=0}^k h(2i)\leq \exp\left\{\displaystyle\frac{2|1-a-\alpha|}{(1-a^2)(1+a)}\right\}.$
\end{enumerate}\end{proposition}
{\bf Proof.} 
\begin{enumerate}

\item[{\rm 1.}] As in the proof of Proposition \ref{pconv}, one gets
$$\displaystyle\prod_{i=0}^k h(2i+1)\leq 
\exp\left\{(1-a-\alpha)(1-a)\sum_{i=0}^{k-1}\frac{a^{2i+1}}{a^{2i+2}(1-a)+\alpha (1-a^{2i+2})}\right\}.$$
First, if  $(1-a-\alpha)a<0$ then 
$$(1-a-\alpha)a(1-a)\sum_{i=0}^{k-1}\frac{a^{2i}}{a^{2i+2}(1-a)+\alpha (1-a^{2i+2})}<0$$ and the result is trivial. 

If $(1-a-\alpha)a\geq 0$ and $a>0$ we have that $1-a\geq\alpha$ and hence
$$a^{2i+2}(1-a)+\alpha (1-a^{2i+2})\geq a^{2i+2}\alpha+\alpha (1-a^{2i+2})=\alpha >0 .$$
Therefore,
\begin{eqnarray*}(1-a-\alpha)a(1-a)\sum_{i=0}^{k-1}\frac{a^{2i}}{a^{2i+2}(1-a)+\alpha (1-a^{2i+2})}\leq\\
(1-a-\alpha)a(1-a)\sum_{i=0}^{k-1}\frac{a^{2i}}{\alpha}<\\
\frac{(1-a-\alpha)a(1-a)}{\alpha}\sum_{i=0}^{\infty}a^{2i}=\frac{(1-a-\alpha)a(1-a)}{\alpha (1-a^2)}<\\\frac{2(1-a-\alpha)a}{1-a^2}=\frac{2|(1-a-\alpha)a|}{1-a^2}.\end{eqnarray*}

On the other hand, if $(1-a-\alpha)a\geq 0$ and $a<0$ then $1-a\leq\alpha$ and this implies
$$a^{2i+2}(1-a)+\alpha (1-a^{2i+2})\geq a^{2i+2}(1-a)+(1-a) (1-a^{2i+2})=(1-a) >0 .$$
Thus,
\begin{eqnarray*}(1-a-\alpha)a(1-a)\sum_{i=0}^{k-1}\frac{a^{2i}}{a^{2i+2}(1-a)+\alpha (1-a^{2i+2})}\leq\\
(1-a-\alpha)a(1-a)\sum_{i=0}^{k-1}\frac{a^{2i}}{1-a}<
(1-a-\alpha)a\sum_{i=0}^{\infty}a^{2i}=\\\frac{|(1-a-\alpha)a|}{1-a^2}\leq \frac{2|(1-a-\alpha)a|}{1-a^2}.\end{eqnarray*}

\item[{\rm 2.}] The same argument  used in the proof of Proposition  \ref{pconv} shows
$$\displaystyle\prod_{i=0}^k h(2i)\leq 
\exp\left\{(1-a-\alpha)(1-a)\sum_{i=0}^{k-1}\frac{a^{2i}}{a^{2i+1}(1-a)+\alpha (1-a^{2i+1})}\right\}.$$

Again, if $1-a-\alpha<0$ then 
$$(1-a-\alpha)(1-a)\sum_{i=0}^{k-1}\frac{a^{2i}}{a^{2i+1}(1-a)+\alpha (1-a^{2i+1})}<0$$
and the inequality of the statement is straightforward.

When we suppose that  $1-a-\alpha\geq 0$ and $a>0$ then, as we did before,
$$a^{2i+1}(1-a)+\alpha (1-a^{2i+1})\geq \alpha >0 $$
and, hence,
\begin{eqnarray*}(1-a-\alpha)(1-a)\sum_{i=0}^{k-1}\frac{a^{2i}}{a^{2i+1}(1-a)+\alpha (1-a^{2i+1})}\leq\\
(1-a-\alpha)(1-a)\sum_{i=0}^{k-1}\frac{a^{2i}}{\alpha}<\frac{(1-a-\alpha)(1-a)}{\alpha (1-a^2)}\leq\\
\frac{(1-a-\alpha)}{ 1-a^2}\leq \frac{(1-a-\alpha)2}{ (1+a)(1-a^2)}.\end{eqnarray*}

Finally, if $1-a-\alpha\geq 0$ and $a<0$, then $a^{2i+1}>a$, and therefore
$$a^{2i+1}(1-a)+\alpha (1-a^{2i+1})=a^{2i+1}(1-a-\alpha)+\alpha >a(1-a-\alpha)+\alpha=(1-a)(a+\alpha).$$
Since $\alpha >(1-a)/2$, we get that $a+\alpha>(1+a)/2$ which leads us to
$$a^{2i+1}(1-a)+\alpha (1-a^{2i+1})>(1-a^2)/2>0.$$
Therefore,
\begin{eqnarray*}(1-a-\alpha)(1-a)\sum_{i=0}^{k-1}\frac{a^{2i}}{a^{2i+1}(1-a)+\alpha (1-a^{2i+1})}\leq\\
(1-a-\alpha)(1-a)\sum_{i=0}^{k-1}\frac{a^{2i}}{(1-a^2)/2}<\frac{(1-a-\alpha)(1-a)}{ (1-a^2)(1-a^2)/2}=
\frac{(1-a-\alpha)2}{(1+a)(1-a^2)}.\end{eqnarray*}
\end{enumerate}\qed\\

\begin{proposition}
\label{pn2} If $0<|a|<1$ and $\alpha >(1-a)/2$, then
\begin{enumerate}
\item[{\rm 1.}] $\displaystyle\prod_{i=0}^k h(2i+1)\geq \min\{1,(1-a)/\alpha\}\exp\left\{-\displaystyle\frac{2|1-a-\alpha|}{(1-a^2)(1+a)}\right\}.$
\item[{\rm 2.}] $\displaystyle\prod_{i=0}^k h(2i)\geq \min\{1,(1-a)/\alpha\}\exp\left\{-\displaystyle\frac{2|(1-a-\alpha)a|}{1-a^2}\right\}.$
\end{enumerate}\end{proposition}
{\bf Proof.} 
First, notice that we have
$$\displaystyle\prod_{i=0}^k h(2i+1)\displaystyle\prod_{i=0}^k h(2i)=\frac{1-a}{a^{2k+2}(1-a)+\alpha (1-a^{2k+2})}.$$
If $\alpha <1-a$ then 
$$\frac{1-a}{a^{2k+2}(1-a)+\alpha (1-a^{2k+2})}>\frac{1-a}{1-a}=1, $$
and when $\alpha \geq 1-a$ we get
$$\frac{1-a}{a^{2k+2}(1-a)+\alpha (1-a^{2k+2})}>\frac{1-a}{\alpha}.$$
Therefore, we get the inequality
$$\displaystyle\prod_{i=0}^k h(2i+1)\displaystyle\prod_{i=0}^k h(2i)\geq\min\{1,(1-a)/\alpha\}$$
and the result claimed follows at once from Proposition \ref{pn1}.\qed\\

\medskip

\begin{theorem} If $0<|a|<1$ then every nonzero periodic solution of (1.1) is stable.
\end{theorem}
{\bf Proof.} 
As mentioned above, if $0<|a|<1$ then every nonzero periodic solution $\{x_n\}$ of (1.1) is given by $x_{2k-1}=p,x_{2k}=q$ for all $k\geq 0,$ where $bpq=1-a$. Let us, then, fix $p,q$ such that $bpq=1-a.$

Since the mapping $f(x_{-1},x_0)=1-a-bx_0x_{-1}$ is continuous, we may find $\delta _1>0$ such that
$|1-a-bx_0x_{-1}|<(1-a)/2$ whenever $||(x_{-1},x_0)-(p,q)||_\infty<\delta_1$.

Let $\{x_n\}$ be the solution of (1.1) obtained for some initial conditions $(x_{-1},x_0)$ verifying  $||(x_{-1},x_0)-(p,q)||_\infty<\delta_1$. According to Propositions \ref{pn1} and \ref{pn2}, there exist continuous functions $f_i(x_{-1},x_0),\ 1\leq i\leq 4$ such that $f_i(p,q)=1$ for $1\leq i\leq 4$ and
\begin{eqnarray*}
& & f_1(x_{-1},x_0)\leq \displaystyle\prod_{i=0}^k h(2i+1)\leq  f_2(x_{-1},x_0),\\
& & f_3(x_{-1},x_0)\leq \displaystyle\prod_{i=0}^k h(2i)\leq  f_4(x_{-1},x_0),
\end{eqnarray*}
for all $k\geq 0$.
For every $\epsilon>0$, we can therefore find $\delta _2>0$ such that $||(x_{-1},x_0)-(p,q)||_\infty<\delta_2$ implies
$$|f_i(x_{-1},x_0)-1|<\frac{\epsilon}{2M}\, ,\quad 1\leq i\leq 4,$$
where $M=\max\{|p+\epsilon|,|p-\epsilon|,|q+\epsilon|,|q-\epsilon|\}$. This clearly implies that, for every $k\geq 0$,
$$\left|\displaystyle\prod_{i=0}^k h(2i+1)-1\right|<\frac{\epsilon}{2M}\,\; ,\quad \left|\displaystyle\prod_{i=0}^k h(2i)-1\right|<\frac{\epsilon}{2M}.$$

Since, by Theorem \ref{t1},
$$x_{2k}=x_0\prod_{i=0}^k h(2i+1),$$ it follows that 
$$|x_{2k}-x_0|
=|x_0|\ \left|\displaystyle\prod_{i=0}^k h(2i+1)-1\right|\leq |x_0|\, \frac{\epsilon}{2M}.$$
Now, if we choose $\delta =\min\{\delta _1,\delta _2,\epsilon/2\}$ and $||(x_{-1},x_0)-(p,q)||_\infty<\delta,$ then we have 
$q-\epsilon< x_0< q+\epsilon$ which implies $|x_0|<M,$ and then
$$|x_{2k}-q|\leq |x_{2k}-x_0|+|x_0-q|\leq |x_0|\frac{\epsilon}{2M}+\frac{\epsilon}{2}<\epsilon.$$

The same argument applied to the subsequence $\{x_{2k+1}\}$ completes the proof.\qed\\

%%%
%%%
%%%
\section{Conclusions and open problems}
%%%
%%%
%%%

We described completely the dynamics and stability properties of Equation (\ref{1}) for all values of the real coefficients $a,b$. This was possible because of the relation of this equation with the M\"obius recurrence  (\ref{2}). We list some open problems related to Eq. (\ref{1}).
\begin{enumerate}
\item Analyze the behaviour of the solutions of (\ref{1}) when the coefficients and the initial conditions are complex. For a recent work on the dynamics of  M\"obius transformations with complex coefficients, see \cite{cgm}.
\item In the case $|a|<1$, try to determine the actual value of the $2$-periodic point $(p,q)$ to which a regular solution of (\ref{1}) converges for a given initial condition $(x_{-1},x_0)$. In the light of Theorem \ref{t1},
this is equivalent to find the value of the infinite products $\prod_{i=0}^{\infty}h(2i)$, $\prod_{i=0}^{\infty}h(2i+1).$ 
\end{enumerate}

%%%
%%%
%%%
\section*{Acknowledgements}
%%%
%%%
%%%

This research was supported in part by the Spanish Ministry of Science and Innovation  (formerly {\it Ministry of Science and Education})   and FEDER,  grant
MTM2007-60679.


\begin{thebibliography}{10}
\bibitem{ag} R. P. Agarwal,  {\it Difference Equations and Inequalities. Theory, Methods, and Applications,} Second edition, Monographs and Textbooks in Pure and Applied Mathematics, {\bf 228}, Marcel Dekker, Inc., New York, 2000.
\bibitem{al} M. Aloqeili,  
Dynamics of a rational difference equation,
Appl. Math. Comput.  {\bf 176} (2006), 768--774.
\bibitem{am2}A. Andruch-Sobilo and M. Migda, Further properties of the rational recursive sequence $x\sb {n+1}={ax\sb {n-1}\over b+cx\sb nx\sb {n-1}}$, Opuscula Math.  {\bf 26} (2006),  387--394. 
\bibitem{am} A. Andruch-Sobilo and M. Migda, 
On the rational recursive sequence {$x\sb
        {n+1}=\frac{ax\sb {n-1}}{b+cx\sb n x\sb {n-1}}$}, 
        Preprint.
        \bibitem{cgm} A. Cima, A. Gasull, and V. Ma\~{n}osa,  Dynamics of some rational discrete dynamical systems via invariants,  Internat. J. Bifur. Chaos Appl. Sci. Engrg.  {\bf 16}  (2006), 631--645.
\bibitem{ci1} C. {\c{C}}inar, 
On the positive solutions of the difference equation {$x\sb
        {n+1}=(x\sb {n-1})/(1+x\sb nx\sb {n-1})$}, 
        Appl. Math. Comput.  {\bf 150} (2004), 21--24.
\bibitem{ci2} C. {\c{C}}inar, 
On the positive solutions of the difference equation {$x\sb
              {n+1}=\frac{ax\sb {n-1}}{1+bx\sb n x\sb {n-1}}$}, 
              Appl. Math. Comput.  {\bf 156} (2004), 587--590.
\bibitem{ci3} C. {\c{C}}inar, 
On the solutions of the difference equation {$x\sb{n+1}=\frac{x\sb {n-1}}{-1+ax\sb nx\sb {n-1}}$}, 
Appl. Math. Comput.  {\bf 158} (2004), 793--797.
\bibitem{ci4} C. {\c{C}}inar, 
 On the positive solutions of the difference equation {$x\sb
              {n+1}=\frac{x\sb {n-1}}{1+ax\sb n x\sb {n-1}}$}, 
       Appl. Math. Comput.  {\bf 158} (2004), 809--812.
\bibitem{ci5} C. {\c{C}}inar, 
On the difference equation {$x\sb {n+1}=\frac{x\sb{n-1}}{-1+x\sb nx\sb {n-1}}$},
 Appl. Math. Comput.  {\bf 158} (2004), 813--816.
\bibitem{cull} P. Cull, M. Flahive and R. Robson,  
{\it Difference Equations. From Rabbits to Chaos,}   Springer, New York, 2005.
\bibitem{el} S. Elaydi,  
{\it An Introduction to Difference Equations,} Third Edition,   Springer, New York, 2005.
\bibitem{kl} M.~R.~S.~Kulenovi\'c and G. Ladas,   {\it Dynamics of Second Order Rational Difference Equations,}   Chapman \& Hall/CRC, Boca Raton, 2002.
\bibitem{st} S. Stevi{\'c}, 
More on a rational recurrence relation,
Appl. Math. E-Notes  {\bf 4} (2004), 80--85.
\end{thebibliography}
\end{document}